**Werner DePauli-Schimanovich**
**Institute for Information Science, Dept. DB&AI, TU-Vienna, Austria**
**Werner.DePauli@gmail.com**

## Naïve Axiomatic Mengenlehre for Experiments (Paper), K50-Set7b

### (0) Abstract

The main goal of "Naïve Axiomatic Mengenlehre" (NAM) is to find a more or less adequately explicit criterion that precisely formalizes the intuitive notion of a "normal set". NAM is mainly a construction procedure for building several formal systems NAMix, each of which can turn out to be an adequate codification of the contentual naïve set theory. ("i" is a natural number which enumerates the used "normality" condition, and "x" is a letter which points to the variants of the used axioms.) Parallel to NAM, the Naïve Axiomatic Class Theory NACT is constructed as a system of systems too.[1]

### (1) The System NAM of Systems NAMix

In NAM any arbitrary formula can be used as a starting point for the investigation of the "normality" of sets. Once a predicate has been found that is equivalent to the predicate atom "Normal" (and which does not falsify NAM of course) then the concept of a normal set is characterized correctly and uniquely.[2] The tack taken here is similar to attemps of the '30s of the past century to find an adequate explication (or surrogate) for the concept of computability. In NAM one can experiment and study the results and effects of different normality conditions until we have found the most appropriate one. With such a condition mathematicians would finally have a decision criterion for generating sets of the naïve set theory.

In the following we want to suggest a series of criteria as conditions for the characterization of normal sets. If we consider also combinations of them, there will be about 100 criteria. One of them will fit almost certainly. Then mathematicians and logicians can prove formally that a given set is normal. Following this line naïve set theory turns out to be completely free of antinomies or paradoxes once we interpret the comprehension schema correctly. Only the normal sets are predicate extensions! The abnormal ones are e.g. identical to the universal set, or to the empty set, or to some external constant "@", just as we like!

Now let us start with the axiomatic construction of all these systems of NAM:
The only objects in NAM are the sets we want to denote by small Latin letters x, y, z, . . . . . .
We also use the (usual unary) prime predicate "Normal". The set operator (S-O) produces always only a set (as it is usual in naïve set theory):
{x: A(x)}  is a SET !

---

[1] In another forthcoming paper, the most important systems of NACT are explained.
[2] Of course it is also possible to construct a hierarchy of characterizations with increasing strength, excluding cases which are too strongly counter-intuitive.

(And that holds for all A in most systems of NAM.)
But we have to distinguish the set operator from the class operator (C-O):
{|x: A(x)|}  is a CLASS !
(Keep in mind that classes are denoted by capital Latin letters X, Y, Z, . . . . . , while the small variables are restricted to sets. But we do not use classes in NAM; and we mentioned this only to avoid confusion! (CoS) for classes cannot produce a contradiction because the pathological predicates A only yield proper-classes and no sets. Therefore we want to call this scheme Church Scheme (CS) for classes, in distinction to (CoS) for sets!

In NAM the general (CoS) is used only in a restricted form. Let "?" be an arbitrary junctor. Then the most important junctive restriction of (CoS) is:
(Ju-CoS) := forall wff A: forall y: [y in {x: A(x)} <==> A(y) "?" Normal] .
"Normal" can also be a structural property of the wff A, e.g.: "Normal:= {x:A(x)} is a CHnP-PE" (Closed Hereditary-non-Patho) yields NAM* [See also K50-Set4]; "Normal := Stratified" yields NF; "Normal := (A is a ZF-axiom)" yields ZF.

## (2) The most important Comprehension Scheme (raBaDi-CoS) of NAM

Generally Normal is a unary predicate of FOL with the basic set of (CoS) substituted into its free variable. The most frequent junctor "?" is "or non" added directly after A on the right side of the equivalence. That gives the most important instance of (Ju-CoS):
(raBaDi-CoS) := (y in {x: A(x)} <===> [A(y) or non-Normal({x: A(x)}]).

We want to call this axiom the "rightside-abnormal-Basicset-Disjunction" Comprehension-Schema, shortly: (raBaDi-CoS).[3] With (raBaDi-CoS) the normal sets are predicate extensions as usual. But the abnormal sets blow up to a set which is "extensional-equal" to the universal set "us := {x: Verum}". They behave somehow like "Urelemente" which are at the same time sets. (Therefore let us call them set-like urelements.) We have 3 possibilities to handle them:

    (A) We identify all pathological sets with the universal set, and the usual axiom of extensionality (EE) works unrestricted as before:
x = y <===> forall z: (z in x <=> z in y) .
But because of substitutivity the universal set has to be abnormal. (But this merely philosophical disadvantage will be considered as a benefit by many readers.)
    (B) The second possibility would be to restrict the extensionality axiom, similarly to what we do in ZFU (= ZF with Urelemente), i.e. to connect the axiom with the condition "Normal(x) and Normal(y)". Then one cannot deduce the equality of the Russell set "ru" with the universal set "us", and "us" can be considered as normal from now on. But this solution complicates the machinery in an unnecessary way.
    (C) The third possibility is a restriction of substitutivity. But since this axiom is part of pure logic, this possibility is certainly the least desirable way out (to solve this problem). Therefore we decide to use at the beginning possibility (A).

---

[3] The author has investigated several restrictions of (CoS) but this one is the best! Of course, one can add Normal also to the left side of the equivalence, or put it as a condition on the entire formula (CoS).

One of our wishes is to find out if we can use complements too, i.e.
(KN-EA) Normal(x) ==> Normal (ko(x)) should be valid.

Therefore also "us" should be normal what is only possible if we use the philosophy of (B). In this case we shall use instead of (EE) the extensionality definition (DefEx): x extensional y [or: x is extensional equal to y]
x ≈ y : = forall z : (z in x <=> z in y).
The equality will be defined by (DefEq):
x = y: = x ≈ y and ([Normal(x) and (Normal(y)] or [non-Normal(x) and non-Normal(y)])

Only too different types of sets can never be equal. Of course (DefEq) can also be formulated as axiom (EqA) if equality "=" is used as a basic symbol of logic.

(raBaDi-CoS) is so important that we want to investigate it in more detail. "Antinomies arise by construction of sets which are extensionally equal to the universe.", said Abraham Fraenkel. And Jonny von Neumann concluded in the same spirit: "If they are not already equal to the universe, the antinomial sets should be at least of equal power!" (In the new language of class theory they are of course the proper classes and the universal class, respectively.)
 (raBaDi-CoS) does more than Jonny wanted. Using philosophy (A) it unifies all abnormal sets with the universal set us := {x: Verum}, what will certainly warm the hearts and enjoy the ZF people. But only for a moment, because (in many systems of NAM) all derivatives of the universal set are normal again, e.g. {x: x =/= us}.

(raBaDi-CoS) is very strong, and together with the axiom of extensionality (EE), a suitable axiom of choice (AC), the set operator (S-O), and some Fundamental- and Eventual-Axioms, it forms a formal system NAM% which is (together with a suitable Normal Condition) similar to naïve set theory and the mathematical praxis (as well as NAM*).

## (3) Other Comprehension Schemata of NAM
It is not an accident that (raBaDi-CoS) forms a system of equal strength as NAM%. There exist also other junctive-restricted comprehension schemes of construction series (Ju-CoS) which achieve something similar. For instance, (rinoBaCo-CoS), the "rightside-normal-Basicset-Conjunction" comprehension scheme
y in {x: A(x)} <===> A(y) & Normal({x: A(x)}) .

In the same way as (raBaDi-CoS) inflates the pathological sets to the universal set, (rinoBaCo-CoS) collapses them to the empty set-like urelements. If we once again choose the philosophy (A), it turns out that the empty set 0 := {x: Falsum} becomes abnormal, which, to be sure, is unfamiliar but in no way inconsistent. Already the singleton {0} of 0 is again normal (and also everything construed upon it). In order to incorporate the empty set into our consideration of the normal sets, we would have to add into the set operator and the axioms the term "or y = 0" in each case where the condition Normal(y) appears.

(rinoBaCo-CoS) & (EE) & (AC) & (S-O) and the Fundamental- and Eventual-Axioms form again together with a suitable Normal-Condition a similar system as NAM%, which we want to call NAM&. Its logical machinery works exactly parallel to NAM%.

From both mentioned systems NAM% and NAM& we can imagine that it would be enough to select an arbitrary set and let it become abnormal to avoid the antinomies. To arrange this we use this time around an implicative restriction (I-CoS) of (CoS), e.g. the "normal-Basicset-Implication" Comprehension-Scheme (noBI-CoS):
Normal({x: A(x)}) ===> [y in {x: A(x)} <=> A(y)] .

By specification of the antecedence with a suitable Normal-condition, all systems like NAM% and NAM& which are characterized uniquely by a restricted (CoS) like (raBaDi-CoS) and (rinoBaCo-CoS) can be constructed. Let us choose e.g. the Normal Condition "not-Equal-@-Implies-normal-Basicset" (nE@I-noB):
{x: A(x)} =/= @ ===> Normal({x: A(x)}),
where "@" is any arbitrary set (which we want to become abnormal in the antinomic case, e.g.: ru, on, us, 0, etc). (noBI-CoS) again together with the other above mentioned axioms yields also a system, we want to call NAM§. (We can add also the re-implications (CoS-I-noB) and (noB-I-nE@) and study these systems. But this often turns out to be too strong.)

We can also adjoin the "commercial at" "@" as a primitive constant to the language of logic with the effect that then the non-pathological predicate extensions can always form a normal set and the proper set universe is free from antinomies. Pathological sets are then "aliens" which have nothing to do with set theory any longer.

### (4) Basic Axioms of NAM
The Naive Axiomatic Mengenlehre NAM should consist in general of 4 Basic Axioms (BAi), 4 Fundamental Axioms (FAj), some Eventuality Axioms (EAk), and a reasonable Normality Condition (NCl) (selected from the group of normality conditions NC).

The 4 Basic Axioms are:
(BA1) := (S-O). The set operator {x: A(x)} produces a set for each wff A .
(BA2): one of the "normalized" Comprehension Schemata (Ju-CoS) or (I-CoS), e.g.:
(BA2a) := (raBaDi-CoS), or
(BA2b) := (rinoBaCo-CoS), or
(BA2c) := (noBI-CoS), or
(BA2d) := what you like.
In the following we want to concentrate our attention mainly on (raBaDi-CoS).

To the set operator and the scheme we have to add:
(BA3) the Axiom of Extensionality (E), mainly as
(BA3a) in the usual formulation of "Element-Equality" (EE):

forall z: (z in x <=> z in y) ===> x = y,[4]

(BA4) an appropriate formulation of the Axiom of Choice (AC), e.g.:
(BA4a) the Choice-Function Axiom (ACF):
exist f: Function(f) & forall y (=/= empty) Subset x ==> f(y) in y . Or:
(BA4b) the Choice-Set Axiom (ACS):
forall y, z in x: (y =/= empty & y intersection z = empty) ===>
exist u: forall v in x: exist w: u intersection v = {w} . Or:
(BA4c) the Ordinals-Universalset Axiom (On-us), e.g. in the form: On ~ us .
Here "~" is a symbol for equal power, i.e. the existence of a bijection between the 2 sets (or classes). On := {x: On(x)} is the class or set of all ordinal-numbers, and us: = {x: Verum} is the universal set.
(BA4d) some other suitable axiom of choice (e.g. the one restricted to small sets), or maybe none.

## (5) Foundamental Axioms of NAM

Furthermore in NAM the following 4 Fundamental Axioms (FAj) should be valid. They express the standard assumptions for normal sets.[5]
(FA1) the Omega Axiom (Om):
Normal(omega),
(FA2) the Power Axiom (Po):
Normal(x) ===> Normal(p(x)),
where p(x) := {y: y Subset x} is the power-set of x.[6]
Note that not every sub-class Y of x (constructed with the wff A) must also be a subset y of x.
E.g.:Ru =< us\{us} & ru = us =/< us\{us} if you use (raBaDi-CoS).
(FA3) the Union Axiom (Un):
Normal(x) & forall y in x: Normal(y) ===> Normal(union(x)),
where union(x) := {y: exist z: y in z & z in x}.
(FA4) shall be a suitable Axiom of Replacement (or Image-set Axiom in the language of class theory) (Im), that the image of a function with a normal domain is also normal. We distinguish different degrees of strength of (Im).[7]

Let us first list the first 3 versions of the image-set axioms:
    (FA4alfa) the "Simple" Image-set Axiom (S-Im):

---

[4] Sometimes (as e.g. in the GNAM-systems considered later) we use also (BA3b), i.e. the "Normal-Elements-Equality" (NEE): forall z: [Normal(z) ==> (z in x <=> z in y)] ===> x = y.
[5] The predicate "Normal(x)" cannot become totally "Falsum" because the Fundamental Axioms and some Eventuality Axioms contribute implicitly to the Normal-Atom.
[6] p(x) denotes the power set (of a set or class), while P(X) denotes the power class (of all subsets of a given class), and PC(X) should denote the 2nd order power class of a class. Sometimes it may be more suitable to define p(x) := {y : y Subset x & Normal(y)}, e.g. in the GNAM-systems.
[7] A is here (and also in the following axioms of (FA4)) a Functional-Formula, i.e. a formula with the functional-property (similar to those used in ZF). As abbreviation for this we write F-F(A).

Normal(x) & Functional-Formula(A) & [exist u in x: exist v: A(u, v)] &
y = {z: [exist w: w in x & A(w, z)]} ======> Normal(y).
  (FA4beta) the "Normal-or-Null-Elements" Image-set Axiom (NoNE-Im):
Normal(x) & Functional-Formula(A) & [exist u in x: exist v: A(u, v)] &
y = {z: [exist w: w in x & A(w, z)] & [Normal(z) or z = 0]} ======> Normal(y).
  (FA4gamma) the "Normal-Elements" Image-set Axiom (NE-Im):
Normal(x) & y = {z: Functional-Formula(F) & [exist u: u in x & F(u, z)] & Normal(z)}
=======> Normal(y).

This partial system of NAM consisting of (raBaDi-CoS) together with (EE), (AC) and the 4 Foundamental Axioms [without any Eventuality Axioms or Normal Conditions] can be consideres somehow as parallel to NBG class theory : a great part of the normal sets are the sets of ZFC and the abnormal sets the proper classes. We want to call this system therefore NAM-ZF.

## (6) Eventuality Axioms of NAM

After presenting the Basic Axioms (BAi) and the Fundamental Axioms (FAj), we want to formulate now the Eventuality Axioms (EAk)[8]. For some combinations of basic axioms and fundamental axioms (and maybe the used Normal Conditions) it may be clever and wise to add some of these axioms too. E.g.:

  (EA1) the "Singleton-Arbitrary-Normal" axiom (SAN):
Normal({@}) , where @ := ru, on, us, 0, etc.
  (EA2) the "Equivalent-Formulas-Equal-Sets" axiom (EFES):
[forall x: A(x) <=> B(x)] ===> {y: A(y)} = {z: B(z)} ,
  (EA3) the "not-Normal-Equal-Mighty" axiom (nNEM):
non-Normal(x) & non-Normal(y) ===> x ~ y ,
  (EA4) the "Elements-Normal-or-Arbitrary" axiom (ENoA):
Normal(x) ====> forall y: [y in x ==> Normal(y) or y = @] .

In addition to these 4 main axioms we can use also:
  (EA5) the "Elements-Normal" Eventuality Axiom (ElN-EA) instead of (EA4):
Normal(x) ===> forall y: [y in x => Normal(y)] ,
  (EA6) the "Komplement-Normal" Eventuality Axiom (KN-EA), which postulates the normality of the complement of a normal set:
Normal(x) ===> Normal(ko(x)).
  (EA7) the "Normal-Set-exclusive-or-Komplement" axiom (NSxorK):
Normal(x) <=/=> Normal(ko(x)).
  (EA8) the"Normal-Set-or-Komplement" axiom (NSoK) axiom:
Normal(x) or Normal(ko(x)),
this axiom plays an important role, but it contradicts (ElN-EA).
  (EA9) the "Normal-Set-or-Supplement" (NSoS) axiom is also conceivable:
Normal({x: A(x)}) or Normal({x: Supplement-A(x)}),

---

[8] "Eventuality" means case-based in this article.

where Supplement-A is that formula, where the Epsilon (i.e. the element sign) is replaced everywhere by its negation.

(EA10) Also other similar conditions on wffs A would make sense to investigate.

As already mentioned above, the systems of NAM are constructed in general as a combination of 4 groups of axioms: basic, fundamental, eventuality axioms and a normal condition, where we sometimes have to choose the most suitable variant of an axiom. But before we start to consider the first normal conditions (which are in fact only the objects for our experiments to find out the most appropriate formal explication of the notion "Normal") the reader should ask himself: what follows already directly from the axioms of NAM without the normal conditions. We want to call these "conditionless" kernel systems NAM0x, where the "x" symbolizes the construction series (i.e. the variants of the (BA2) axioms). "a" points e.g. for the use of (raBaDi-CoS), "b" for that of (rinoBaCo-CoS), "c" of (noBI-CoS), etc. (NAM0a, b, c are more or less other names for NAM%, &, §.)

## (7) The first Normal Conditions and the first Systems of NAM

After this short interruption let us go on with the explanation of the last group: the Normal Conditions (NCl). We want to find the best of these conditions to select. But until we know which condition is the best, we will need some time for experiments. The first two conditions are:

(NC1) the "not-Function-Domain-Komplement" Condition (nFDK-Cond):
non-exist f: (Function(f) & f[x] = ko(x)) ===> Normal(x) or x = us.
f[x] is the image of the domain x under f (and not the value f(x)), and ko(x) is the Complement {y: y non-in x}. By didactical reasons we can add the term "exist g: (Function(g) & g[ko(x)] = x)" in the antecedent, but this term is in fact redundant.

(NC2) the "Smaller-Power-than-Komplement" Condition (SPK-Cond):
card(x) < card(ko(x)) [& x ~/~ ko(x)] ===> Normal(x) ,
where card(x) is the cardinal-number of x, and card(ko(x)) that of the complement. (Again the didactical term in brackets can be omitted.)
(NC1) and (NC2) are essentially only different formulations of one and the same fact.[9]

If we keep the above-mentioned connections in mind, we get first the semi-canonical system NAM1a with (BA2a) := (raBaDi-CoS) plus the other BAi, and (FA4alfa) := (S-Im) plus the other FAj, with condition (NC1). The first 3 EAl are derivable in NAM1a. This system NAM1a will be later the basis for some extensions, e.g. for NAM1aKNoU := NAM1a & (KNoU-EA).

---

[9] The danger inherent in these conditions can be explained by the following example: Let us define Slim := {x: card(x) < card(ko(x))} and HeredSlim := {x: forall y in x*: card(y) < card(ko(y))}. We want to find out if "Slim in Slim" or "HeredSlim in HeredSlim" ? That would imply an antinomy similar to Cantor's. Fortunately the "Smaller-Power-than-Komplement" Condition (SPK-Cond) and "Hull-Smaller-Mightiness" (HSM-Cond) are consequently wrong for Slim and HeredSlim.

In analogy to NAM1a we get the parallel-semicanonical system NAM1b which is also founded on condition (NC1), but together with (BA2b) := (rinoBaCo-CoS) and (FA4beta) := (NE-Im) and the other axioms (like in NAM1a) plus the 3 eventuality axioms. Like NAM1b we get the quasi-semicanonical System NAM1c with (BA2c) := (noBI-CoS), (FA4beta) := (NE-Im), as well as (NC1), plus the other axioms (inclusive the 3 eventuality axioms). (NAM1c does not imply Cantor's Antinomy.) NAM2a, b, und c, are the same systems as NAM1a, b und c, only based on condition (NC2). NAM2c is a better basis for extensions than NAM2a, yielding NAM2cKN := NAM2c & (KN-EA). Thus, we get (with this correspondence of the axioms) all together 6 basic systems (in addition to the 3 kernel systems). Note also that NAM2a without foundational axioms, but together with (EA7) := (NSxorK) is equivalent to NACT#. Together with the 4 foundational axioms it is equivalent to NACT#4. (See K50-Set3b.)

Let us first write down the axioms of the first system NAM1a of the Naïve Axiomatic Mengenlehre in a list. NAM1a has 8 axioms:
(1-1) {x: A(x)} is a set,
(1-2a) [y in {x: A(x)} <==> A(y) or non-Normal({x: A(x)})] ,
(1-3) forall z: (z in x <==> z in y) ===> x = y,
(1-4) a suitable axiom of choice, e.g. On := {|x: On(x)|} ~ us,
(1-5) Normal(omega) ,
(1-6) Normal(x) ==> Normal(p(x)) ,
(1-7) Normal(x) & forall y in x: Normal(y) ==> Normal(union(x)),
(1-8) Normal(x) & Functional-Formula(A) & [exist u in x: exist v: A(u, v)] &
y == {z: [exist w: w in x & A(w, z)] & Normal(z)} ====> Normal(y) ,
(1-9) non-exist f: (Function(f) & f[x] = ko(x)) ===> Normal(x) or x = us.

Axiom (1-9) is (NC1). The leading "1" points at the 1st NC, because we are enumerating the systems synchronously with the Normal Conditions.
The Eventuality Axioms (1-10) to (1-13) are derivable from NAM1a.
(1-10) Normal({us}) ,
(1-11) [forall x: A(x) <=> B(x)] ===> {y: A(y)} = {z: B(z)} ,
(1-12) non-Normal(x) & non-Normal(y) ===> x ~ y ,
(1-13) Normal(x) ==> forall y: [y in x ==> Normal(y) or y = us] ,
NAM1a produces mainly such normal sets which are similar to the sets of ZFC; but also much more than that.

### (8) Complements and restricted Image-set Axioms
If we want to add a complement to NAM1a we would put also the "Komplement-Normal-or-Universal" Eventuality Axiom (EA6') := (KNoU-EA) as axiom with the number 14 onto the list:
(1-14) Normal(x) ==> Normal(ko(x)) or ko(x) = us,

Let us call this system NAM1aKNoU, where the pathological sets in the "middle" are abnormal, while small and large sets (without "us") are normal.[10] But in that case we have to restrict the domain in the Image-set Axioms (1-7alfa, beta, und gamma) to small sets; otherwise it would be possible to deduce the normality of the Russell-set "ru". Therefore (1-8) has to be replaced by stronger restricted axioms (e.g. (FA4delta, eta, fi, psi, chi, jota, or kappa)). Before we can decide which one we want to select let us enumerate the new image-set axioms in a list.

(FA4delta) the "twofold Functional-Formula" Imageset Axiom (2FF-Im):
Normal(x) & y1 ==
== {y: ( [Funktional-Formel(F) & Funktional-Formel(G) & (forall x1, y1, u, v: x1 =/= y1 & G(x1, u) & G(y1, v) => u =/= v)]   ====>   [exist z: z in x & F(z, y) & Normal(y)   <==> non-forall w non-in x: exist z in x: G(z, w)] }      ========>   Normal(y1) .
In fact this imageset axiom is a double formula scheme[11] which expresses the same fact as the following two other formulations:

(FA4eta) the "not-Function-Domain-Komplement" Image-set Axiom (nFDK-Im):
Normal(x) & y1 ==
= {y: Functional-Formula(F) & (non-exist g: Function(g) & g[x] = ko(x)) & [exist z: z in x & F(z, y)] & Normal(y)}   ========> Normal(y1) .
Here g[x] is the image of the domain x of the function g.

(FA4phi) the "Domain-Smaller-Power-as-Komplement" Image-Set-Axiom (DSPaK-Im):
Normal(x) & y1 ==
== {y: Functional-Formula(F) & [card(x) < card(ko(x))] & [exist z: z in x & F(z, y)] & Normal(y)}   ========> Normal(y1) .

The condition slim(x) := [card(x) < card(ko(x))] expresses the fact that x should be "small" in some specific sense. But instead of using this formulation of "small" we can also use another condition, e.g.: Mirimanoff, Founded, Hereditary Founded, or Cantorian, which yields the following axioms. (Note that in NAM it depends of the system chosen whether these conditions are equivalent or not.)

---

[10] For NAM1aKN, the same is valid as for ZFK (respectively ZFCK). See Goldstern.
[11] For unsuitable F or G, y1 will yield the universal set. (Therefore it would be good to use either (KNoU-EA) or NAM2c as basis.) The right-side expression in the parentheses (in the antecedent of the implication in the set-brackets) is only added for didactical reasons, since everybody knows that G is an injection. In the consequent (i.e. right side of the implication) in the set-brackets, the contra-valent (respectively the equivalence with the negation) is valid! Therefore not both conditions in the 2nd brackets [= square brackets] can simulaneously have the same truth-value. The term on the right side of the equivalence says that G is no surjection of x onto its complement. E.g.: Let x := omega. Also if you can find an appropriate F-F(G) (such that the antecedent becomes true) where G can map omega onto its complement, the 2nd condition of the contra-valent always becomes wrong (because that is impossible). Therefore omega is of smaller power than its complement, and therefore the 2nd condition is true: we can produce an image of omega! But if x = ko(omega), then the 2nd condition is true and the 1st one must be wrong: we cannot construct an image of ko(omega).

 (FA4psi) the "Mirimanoff-Domain" Image-set Axiom (MD-Im) with:
Mirimanoff(x) := non-exist f: [Function(f) & f(0) = x & forall n in omega: f(n+1) in f(n)].
 (FA4chi) the "Founded-Domain" Image-set Axiom (FD-Im) with:
Found(x) := exist y in x: x intersection y = empty .
 (FA4jota) the "Hereditary-Founded-Domain" Image-set Axiom (HFD-Im) with:
Heri-Found(x) := forall y in e*(x): exist z in y: z intersection y = empty .
 (FA4kappa) the "Cantorian-Domain" Image-set Axiom (CD-Im) with:
Cantorian(x) := [card(x) < card(p(x))] .

If we drop the appendix "Normal(y)" in the curly braces (= set parentheses) of these new 7 axioms, let us call them (FA4d, e, f, g, h, j, und k). This release from "Normal(y)" will allow us to also generate the derivatives of the universal set and similar sets. (In some systems the normality of the elements is redundant because it follows from other facts.) In analogy to NAM1aKN we get further systems like NAM2aKN, etc, if we replace the image-set axiom (1-8) by one of the 7 new (FA4).

To formulate the other NC's we want to explain what is the meaning of the hull e* of the elementhood relation e := {<x, y>: x in y}). Let us define:
x e* y :<==> exist n (=/=0) in omega: exist f: Funktion(f) & f(0)=x & f(n-1)=y & forall i in n: f(i) in f(i+1) .
e* := {<x, y>: x e* y}, [or the union from n = 0 until omega: e^n] .
 It is known that e* can be defined within NBG and ZFC, and -- as one can see -- also within NAM1a, etc. In those systems where that should not be possible, we have to use the class theory CNAM (constructed over NAM) and define the class relations E, E^n and E* in it.

Since with e and E, "e..." and "E..." are also relations, we can fix their ranges as we are accustomed to be "e...(x)" := {y: exist z in x: <z, y> in e...} and "E...(X)" similarly. In analogy to the above we write for the complements of e^n, E^n, e* and E* the following:
n-e^n, n-E^n, n-e* or n-E*, and the same symbols for the predicates.
e^+ (or epsilon-plus) shall be e* \ id and E^+ shall be E* \ Id.

## (9) Further Normal Conditions
With this preliminary work we can proceed with the definitions of the NC's:
 (NC3) is the "Elements-of-Hull-no-Function-from-Domain-to-Komplement" (EHnFDK-Cond):
forall y in e*(x):  [non-exist f: Function(f) & f[y] = ko(y)] &
      [exist g: Function(g) & g[ko(y)] = y]  ====> Normal(x) .
The 2nd term in the brackets in the antecedent is only of didactical nature and can be dropped.
 (NC4) the "Hull-Smaller-Mightiness" Condition (HSM-Cond):
forall y in e*(x): card(y) < card(ko(y))  ====> Normal(x) .
This condition is more complicated than (NC2), but it does hinder the further construction of "counter-intuitive" sets.

   (NC5) the "no-Bijection-Domain-Komplement" Condition (nBDK-Cond):
non-exist f: (Function(f) & f[x] = ko(x))  or …
… or non-exist g: (Function(g) & g[ko(x)] = x)   ====>   Normal(x) .
   (NC6) the "not-Equal-Mighty-Komplement"[12] Condition (nEMK-Cond):
x ~/~ ko(x)  ===>   Normal(x),
what means that if x is not equal mighty with its complement then it is normal.[13] Note also that NAM6c without foundamental axioms, but together with (EA8) := (NSoK) is equivalent to NACT+. Together with the 4 foundational axioms it is equivalent to NACT+4. (See K50-Set3b.)

Just as before we claim now the validity of the last 2 conditions for all elements of the hull of x. Thus we get:
   (NC7) the "Elements-of-Hull-not-Bijection-Domain-Komplement" Condition (EHnBDK-Cond):
forall y in e*(x): "the antecedens of (NC5) with (x/y)"   ====>   Normal(x) .
   (NC8) the "Elements-of-Hull-not-Equal-Mighty-Komplement"  (EHnEMK-Cond):
forall y in e*(x): y~/~ko(y) ====>   Normal(x) .

There exists of course a whole series of further Normal Conditions that generate also more or less interesting systems. E.g.:
   (NC9) the "not-Elementship-of-ItSelf" Condition (nEIS-Cond)[14]:
x non-in x   ====>    Normal(x) .
   (NC10) the "not-Element-of-It's-Own-Hull" Condition (nEIOH-Cond):

---

[12] We want to use in the context the notion "Equal-Mighty" instead of "equipollent", because it allows us the easier abbreviation. Furthermore, it represents the new language Euro-English. We need not to speak American English in Europe, and British English is too difficult for us.

[13] Be careful with (NC6) := "not-Equal-Mighty-Komplement" (nEMK-Cond) or (NC8) := "Elements-of-Hull-not-Equal-Mighty-Komplement" (EHnEMK-Cond)! Because (raBaDi-CoS) and (rinoBaCo-CoS) force "us" or "0" to be abnormal, we should add the appendix "or x = us or x = 0" to "Normal(x)" at the end of (NC5) to (NC8) and call them (NC5') to (NC8'). Or we have to use (noBI-CoS) and "@ := on". And this implies (because of "not-Normal-Equal-Mighty" (nNEM)) the strong axiom of choice (OnAs) := On ~ us. But it would be much easier if we would use philosophy (B).

[14] In some systems (especially those formed with (NC9) etc) e.g. the following is valid: Let ru := {x: x non-in x} be the Russell-set. The conditional Comprehension Scheme (noBI-CoS) together with its re-implication (CoS-I-noB) makes the implication to an equivalence yielding (noBE-CoS). This scheme allows the deduction of the following facts: ru = ru u {ru} = {x: x non-in x  v  x = ru} . ("u" is here the union, and "v" is the "or".) This is reasonable because the right-hand term of the implication (= the one with the equivalence in (noBI-CoS)) produces a contradiction. This falsifies the antecedent and yields "ru in ru". But at the same time the following is also true: ru = ru \ {ru} = {x: x non-in x & x =/= ru}. You cannot take ru out from itself. The same is valid for the other pathological sets too, like the ordinal-set on, the Mirimanoff-set mi, the Founded-set fu, etc. With ru we have also as many pathological sets as its elements, because we can generate (ru \ {x}) and (ru \ x) for every element x in ru (which are mostly normal sets), and for these new products (ru \ {x}) in (ru \ {x}) and also (ru \ x) in (ru \ x) is valid.

x non-epsilon-plus x  ====>   Normal(x),
"x non-in e^+(x)" means: non-exist n =/= 0 in omega: x in x1 in x2 in . . . in xn = x .
  (NC11) Mirimanoff-Condition: "not-Infinite-Descending-Element-Sequence" (nIDES-Cond):
non-exist f: [Function(f) & f(0) = x & forall n in omega: f(n+1) in f(n)] ===> Normal(x) .
  (NC12) the Foundedness Condition (Found-Cond):
exist y in x: x intersection y = empty  ====>   Normal(x) .

These 4 conditions can again be strengthened by the claim of the validity for all elements of the hull of the former domain:
  (NC13) the "Elements-of-Hull-not-Elements-of-ItSelf" Condition (EHnEIS-Cond):
forall y in e*(x): y non-in y  ====>  Normal(x) .
  (NC14) the "Elements-of-Hull-not-Elements-of-its-Own-Hull" (EHnEOH-Cond):
forall y in e*(x): y non-in e^+(y)  ===>  Normal(x) .
  (NC15)  "Elements-of-Hull-not-Infinite-Descending-Element-Sequence"  (EHnIDES-Cond):
forall y in e*(x): non-exist f: [Function(f) & f(0) = y & forall n in omega: f(n+1) in f(n)] ===> Normal(x) .
  (NC16) "Elements-of-Hull-Founded" Condition (EH.Found-Cond):
forall y in e*(x): exist z in y: z intersection y = empty ===> Normal(x) .

(NCxyz) can be any arbitrary Normal Condition. Because we can also consider systems NAMix where "Normal" is constructed completely differently. (e.g. as already mentioned not depending on the elements y or the constructed set {x: A(x)}, but on the structure of the wff A, as in NF.) To investigate all these conditions is the program of NAM.[15]

After studying NAM on basis of (raBaDi-CoS) the reader will probably agree with the author that the use of philosophy (A) causes a lot of difficulties. Only if we want to investigate and experiment in a ZF-like system as NAM-ZF this philosophy makes sense. But if we want to study systems with a complement (what is our intended goal) it is better to use philosophy (B) with (DefEx) and (EqA). So in fact we should now start a new and construct NAMix under (B). But here is not enough space for this endeavour and therefore it is left as homework for the reader.

## (10) Extensions of NAM

---

[15] Further Normal Conditions can be formulated and constructed with structural constraints on the formula A (like Stratified(A)). But we can also use completely different structural properties like "Supplement" (= Counter-valisation, where all epsilons are replaced by non-epsilons in the expanded formula, and vice versa), Dualisation (where all "and" are replaced by "or" and all general quantifiers by existence-quantifiers, and vice versa), and structural invariants of such constraints. It may also be possible to show with (AC) that some NC's are equivalent, or we can derive for some of the NC's an implication chain already in NAM0a or in other NAMix.

Most systems NAMix described up to now are pure set theory. Of course we can always cover these pure set theory systems with a class theory (similar to NBG), getting the systems CNAMix of CNAM. That is important especially because sets and classes generated with the same wff A often have a different extension. (E.g., in some systems the universal set "us" is different from the universal class UC, which can contain abnormal sets like the Russell-set "ru" too [usually different from "us" in some systems], while "us" may not contain abnormal sets.) For classes over NAM-sets we use capital Latin letters X, Y, Z, … etc as variables and {|x: A(x)|} as class operator. (To mark the C-O we use always the vertical bar together with the set-parantheses.)

But NAM can also be extended into the other direction (i.e. downwards) if we restrict the set-constituting variable "x" in the set operator to normal sets, i.e. use always {x: A(x) & Normal(x)}. For normal sets we introduce as new variables the small German letters german-x, german-y, german-z, … etc (abbreviated as g-x, g-y, g-z, … etc) and the normal-set-operator (NS-O) := {g-x: A(g-x)} or {/x: A(x)/}. (To mark the NS-O we use always the slash together with the set-parantheses.) We want to call this theory then the Germanized Naïve Axiomatic Mengenlehre GNAM and its systems GNAMix.[16] If a majority of the NAMix yields (after their Germanization) a unique same system GNAM! (of course up to equivalence), it would be an indication that this system GNAM! is an adequate representative of the formalization of Naïve Mengenlehre!

With the help of normal sets such systems (generated by normal conditions which allow counter-intuitive sets too) can also be made to work very well. We can define in such systems functions, relations, cardinals and ordinals, equality, subsets, power, etc, as precisely as those used in ZFC. If there is no danger of mixing it up we can drop the German letters and use, as accustomed, the Latin ones. And then everything mathematicians are doing in the Naïve Mengenlehre is as usual; only now they know (because of the experiments and the philosophical justification going hand in hand with them) that what they are doing is consistent, inspite of Gödel and the miscarriage of Hilbert's program.

**References for the Set-Theoretical Articles, K50-Set10**

Bernays, Paul & Fraenkel, Abraham
[1968] Axiomatic Set Theory. Springer Verlag, Heidelberg, New York.
Brunner, Norbert & Felgner, Ulrich
[2002] Gödels Universum der konstruktiblen Mengen. In: [Buldt & al, 2002].
Buldt, Bernd & Köhler, Eckehard & Stöltzner, Michael & Weibel, Peter & Klein, Carsten & DePauli-Schimanovich-Göttig, Werner

---

[16] In this case, ZFC would be for example GNAM0a without Eventuality Axioms. Of course (EFES) and (nNEM) would become redundant (or may have to be rejected), and (ElN-EA) derivable for normal sets from (NE-Im) or (S-Im), because in GNAM the sets in the set operator are restricted to normal ones like the classes in the class operator are restricted to sets. Therefore we could also use GNAM0a as basis for ZFC. In both cases the missing ZF-axioms can be deduced from the Foundational Axioms.

[2002] Kurt Gödel: Wahrheit und Beweisbarkeit, Band 2: Kompendium zu Gödels Werk. oebv&htp, Wien.

Casti, John & DePauli, Werner
[2000] Gödel: A Life of Logic. Perseus Publishing, Cambridge (MA).

Davidson, Donald & Hintikka, Jaakko
[1969] Words and Objections. D. Reidel Publ. Comp., Dordrecht.

DePauli-Schimanovich, Werner (See also: Schimanovich [1971a], [1971b]) and [1981]).
[ca 1990] On Frege's True Way Out. Article K50-Set8b in this book.
[1998] Hegel und die Mengenlehre. Preprint at: http://www.univie.ac.at/bvi/europolis .
[2002a] Naïve Axiomatic Mengenlehre for Experiments. In: Proceedings of the HPLMC (= History and Philosophy of Logic, Mathematics and Computing), Sept. 2002 in Hagenberg/Linz (in honour of Bruno Buchberger).
[2002b] The Notion "Pathology" in Set Theora. In: Abstracts of the Conference HiPhiLoMaC (= History and Philosophy of Logic, Mathematics and Computing), Nov. 2002 in San Sebastian.
[2005a] Kurt Gödel und die Mathematische Logik (EUROPOLIS5). Trauner Verlag, Linz/Austria.
[2005b] Arrow's Paradox ist partial-consistent. In: DePauli-Schimanovich [2005a].
[2006a] A Brief History of Future Set Theory. Article K50-Set4 in this book.
[2006b] The Notion of "Pathology" in Set Theory. Article K50-Set5b in this book.
[2006c] Naïve Axiomatic Class Theory NACT: a Solution for the Antinomies of Naïve "Mengenlehre". Article K50-Set6 in this book.
[2006d] Naïve Axiomatic Mengenlehre for Experiments. Article K50-Set7b in this book.

DePauli-Schimanovich, Werner & Weibel, Peter
[1997] Kurt Gödel: Ein Mathematischer Mythos. Hölder-Pichler-Tempsky Verlag, Wien.

Feferman, Solomon & Dawson, John & Kleene, Stephen & Moore, Gregory & Solovay, Robert & Heijenoort, Jean van,
[1990] Kurt Gödel: Collected Works, Volume II. Oxford University Press, New York & Oxford.

Felgner, Ulrich
[1985] Mengenlehre: Wege Mathematischer Grundlagenforschung. Wissensch. Buchgesellschaft, Darmstadt.
[2002] Zur Geschichte des Mengenbegriffs. In: [Buldt & al, 2002] und Artikel K50-Set3 in this book.

Forster, Thomas
[1995] Set Theory with an Universal Set. Exploring an Untyped Universe. ($2^{nd}$ Edition.) Oxford Science Publ., Clarendon Press, Oxford.

Gödel, Kurt
[1938] The relative consistency of the axiom of choice and of the generalized continuum hypothesis. In: [Feferman & al, 1990].

Goldstern, Martin & Judah, Haim
[1995] The Incompleteness Phenomenon. A. K. Peters Ltd., Wellesley (MA).

Goldstern, Martin
[1998] Set Theory with Complements. http://info.tuwien.ac.at/goldstern/papers/notes/zfpk.pdf

Hallet, Michael
[1984] Cantorian Set Theory and Limitation of Size. Oxford Univ. Press, N.Y. & Oxford.
Halmos, Paul
[1960] Naive Set Theory.    Van Nostrand Company Inc., Princeton (NJ).
Holmes, Randall
[1998] Elementary Set Theory with a Universal Set.
Vol. 10 of the Cahiers du Centre de Logique. Academia-Bruylant, Louvain-la-Neuve (Belgium).
[2002] The inconsistency of double-extension set theory.
http://math.boisestate.edu/~holmes/holmes/doubleextension.ps
Jech, Thomas
[1974] Procedings of the Symposium in Pure Mathematics (1970), Vol. XIII, Part 2, AMS, Provicene R.I. .
Jensen, Ronald Björn
[1969] On the consistency of a slight (?) modification of Quine's New Foundation.
In: [Davidson & Hintikka, 1969].
Köhler, Eckehart & Weibel, Peter & Stöltzner, Michael & Buldt, Bernd & Klein, Carsten & DePauli-Schimanovich-Göttig, Werner
[2002] Kurt Gödel: Wahrheit und Beweisbarkeit, Band 1: Dokumente und historische Analysen. Hölder-Pichler-Tempsky Verlag, Wien.
Kolleritsch, Alfred & Waldorf, Günter
[1971] Manuskripte 33/'71 (Zeitschrift für Literatur und Kunst). Forum Stadtpark, A-8010 Graz, Austria.
Mathias, A.R.D.
[1992] The Ignorance of Bourbaki. In: The Mathematical Intelligencer 14 (No.3).
Quine, Willard Van Orman
[1969] Set Theory and its Logic. Belknap Press of Harvard University Press, Cambridge (MA).
Rubin, Jean & Rubin, Herman
[1978] Equivalents of the Axiom of Choice.    Springer, Heidelberg & New York
Schimanovich, Werner,
[1971a] Extension der Mengenlehre. Dissertation an der Universität Wien.
[1971b] Der Mengenbildungs-Prozess. In: [Kolleritsch & Waldorf, 1971].
[1981] The Formal Explication of the Concept of Antinomy. In: EUROPOLIS5, K44-Lo2a (and Lo2b). Or: Wittgenstein and his Impact on Contemporary Thought, Proceedings of the 2$^{nd}$ International Wittgenstein Symposium (29$^{th}$ of Aug. to 4$^{th}$ of Sept. 1977), Hölder-Pichler-Tempsky, Wien 1978.
Schimanovich-Galidescu, Maria-Elena
[2002] Princeton – Wien, 1946 – 1966. Gödels Briefe an seine Mutter. In: [Köhler & al, 2002].
Scott, Dana,
[1974] Axiomatizing Set Theory. In: Jech [1974] .
Suppes, Patrick,
[1960] Axiomatic Set Theory. D. Van Nostrand Company Inc., Princeton (NJ).

Weibel, Peter & Schimanovich, Werner
[1986] Kurt Gödel: Ein mathematischer Mythos. Film, 80 minutes, copyright ORF (= Austrian Television Network), ORF-shop, Würzburggasse 30, A-1130 Wien.
Wittgenstein, Ludwig,
[1956] Bemerkungen über die Grundlagen der Mathematik / Remarks on the Foundation of Mathematics. The M.I.T. Press, Cambridge MA and London.